\magnification1090
 \def\oui{oui}
 \ifx\textures\oui
 \input CrayolaColors 
 \long\def\rge#1{\Red#1\Black}

 \else 

  
  %
  %

  %
  %

  \catcode`@=12 

 \def\defrefnote#1{\definexref{#1}{{\the\footnotenumber}}{refnotes}}

  %
  %


 \input eplain.tex
\makeatletter
\def\numberedfootnote{%
ÊÊ\global\advance\footnotenumber by 1
ÊÊ\@eplainfootnote{{\number\footnotenumber}}%
}%
\def\makecolumns#1/#2 {\par \begingroup
ÊÊ \@columndepth = #1
ÊÊ \advance\@columndepth by -1
ÊÊ \divide \@columndepth by #2
ÊÊ \advance\@columndepth by 1
ÊÊ \@linestogoincolumn = \@columndepth
ÊÊ \@linestogo = #1
ÊÊ \currentcolumn = 1
ÊÊ \def\@endcolumnactions{%
ÊÊÊÊÊÊ\ifnum \@linestogo<2
ÊÊÊÊÊÊÊÊ \the\crtok \egroup \endgroup \par 
ÊÊÊÊÊÊ\else
ÊÊÊÊÊÊÊÊ \global\advance\@linestogo by -1
ÊÊÊÊÊÊÊÊ \ifnum\@linestogoincolumn<2
ÊÊÊÊÊÊÊÊÊÊÊÊ\global\advance\currentcolumn by 1
ÊÊÊÊÊÊÊÊÊÊÊÊ\global\@linestogoincolumn = \@columndepth
ÊÊÊÊÊÊÊÊÊÊÊÊ\the\crtok
ÊÊÊÊÊÊÊÊ \else
ÊÊÊÊÊÊÊÊÊÊÊÊ&\global\advance\@linestogoincolumn by -1
ÊÊÊÊÊÊÊÊ \fi
ÊÊÊÊÊÊ\fi
ÊÊ }%
ÊÊ \makeactive\^^M
ÊÊ \letreturn \@endcolumnactions
ÊÊ \@columnwidth = \hsize
ÊÊÊÊ \advance\@columnwidth by -\parindent
ÊÊÊÊ \divide\@columnwidth by #2
ÊÊ \penalty\abovecolumnspenalty
ÊÊ \noindent 
ÊÊ \valign\bgroup
ÊÊÊÊ &\hbox to \@columnwidth{\strut \hsize = \@columnwidth ##\hfil}\cr
}%
\makeatother

\lefteqnumbers
   \def\testd{oui}
   \def\choixlat{\ifx\numadroite\testd\righteqnumbers
            \else  \lefteqnumbers\fi}
    \choixlat

\catcode`@=\letter
\def\@eplainfootnote#1{\let\@sf\empty 
  \ifhmode\edef\@sf{\spacefactor\the\spacefactor}\/\fi
  \global\advance\hlfootlabelnumber by 1
  \hlstart@impl{foot}{\hlfootlabel}%
  \hldest@impl{footback}{\hlfootbacklabel}%
  \hbox{$^{(#1)}$}%
  \hlend@impl{foot}%
  \@sf\vfootnote{#1.}%
}%
\catcode`@=\other

  \interfootnoteskip=0pt
  \let\note=\numberedfootnote
  \everyfootnote={\eightpoint\leftskip=5truemm\rightskip5truemm}
  
  \hsize150truemm\vsize 240truemm\hoffset=5truemm
  \def\dimstand{\hsize 150truemm\vsize 240truemm\hoffset=5truemm\voffset=0truemm}
  
  \pretolerance=500\tolerance=1000\brokenpenalty=5000
  \parindent3mm
  
  \countdef\temps=170
  \temps=\time
  \countdef\nminutes=171{\nminutes = \time}
  \countdef\nheure=172
  \def\heure{\begingroup                     
     \temps = \time \divide\temps by 60
     \nheure = \temps                        
     \nminutes = \time
     \multiply\temps by 60
     \advance\nminutes by -\temps            
     \ifnum\nminutes<10 \toks1 = {0}%
     \else\toks1 = {}%
     \fi
     \number\nheure h\the\toks1 \number\nminutes  
  \endgroup}%

  \newcount\chstart
  \chstart=\pageno
 \headline={\ifnum\pageno=\chstart {\hfill} \else {\hss \tenrm --\ \folio\ --\hss}\fi}
  \footline={\hfill}
  \normalbaselines
  \frenchspacing
    \def\dater{\vglue-10mm\rightline{(\the\day/\the\month/\the\year)}}
  \def\dateheure{\vglue-10mm\rightline{(\the\day/\the\month/\the\year,\ \heure)}}

  \newif\ifpagetitre \pagetitretrue
\newtoks\hautpagetitre \hautpagetitre={\hfill}
\newtoks\baspagetitre \baspagetitre={\hfill}
\newtoks\auteurcourant \auteurcourant={\hfill}
\newtoks\titrecourant \titrecourant={\hfill}
\newtoks\hautpagegauche
\newtoks\hautpagedroite
\newtoks\hautpagemilieu
\hautpagemilieu={\tenrm\hfil -- \folio\ -- \hfil}
\hautpagegauche={\ifx\midfolio\oui\the\hautpagemilieu\else\tenrm\folio\hfill\the\auteurcourant\hfill\fi}
\hautpagedroite={\ifx\midfolio\oui\the\hautpagemilieu\else\hfill\the\titrecourant\hfill\tenrm\folio\fi}
\newtoks\baspagegauche \baspagegauche={\hfil}
\newtoks\baspagedroite \baspagedroite={\hfil}
\headline={\ifpagetitre\the\hautpagetitre
\else\ifodd\pageno\the\hautpagedroite\else\the\hautpagegauche\fi\fi }
\footline={\ifpagetitre\the\baspagetitre
\else\ifodd\pageno\the\baspagedroite
\else\the\baspagegauche\fi\fi \global\pagetitrefalse}

\def\pageblanche{\vfill\eject\pagetitretrue
\null\vfill\eject
\pagetitretrue
}
\def\chgtpage{\ifodd\pageno \else
\pageblanche \fi \pagetitretrue\titreun=0\footnotenumber=0}

\def\chgtpageincrtitreun{\ifodd\pageno \else
\pageblanche \fi \pagetitretrue\footnotenumber=0}

\def\majnombres{\ifodd\pageno \else
\pageblanche \fi \pagetitretrue\hautpoly\titreun=0\footnotenumber=0}

\def\hautspages#1#2{\auteurcourant={\ninepcap#1}\titrecourant={\nineit#2}}

\ifnum\chstart=\pageno \pagetitretrue\fi
  

  \def\PAR{\par}
  
  \def\leftnote#1{\vadjust{\setbox1=\vtop{\hsize 20mm\parindent=0pt\eightpoint
  \baselineskip=9pt\rightskip=4mm plus 4mm\vskip-4.7mm#1}\hbox{\kern-2cm\smash{\box1}}}}

  
  \def\raggedcenter{\leftskip=20pt plus 10em  
       \rightskip=\leftskip 
        \parfillskip=0pt 
         \spaceskip=.3333em \xspaceskip=.5em 
          \pretolerance=9999 \tolerance=9999
           \hyphenpenalty=9999 \exhyphenpenalty=9999 }
           
  \def\titrecentre#1{{\parindent0mm\raggedcenter
       \spaceskip=.6em plus .2em minus .2em\xspaceskip=.6em plus .2em minus .2em
        \tit#1\par}}
        


  \def\oui{oui}
  
   \def\fontetitreun{\ifx\paradouze\oui\douzepts\gpdouze\twelvebf\textfont1=\twelveib\else
\quatorzepts\gpquatorze\fourteenbf\fi}

\def\fontetitreunl{\douzepts\textfont1=\twelveib\scriptfont1=\tenib\fourteenti}
 
 \def\fontetitredeux{\textfont1=\eleveni\ifx\paradouze\oui\onzepts\scriptfont1=\ninei\elevenit\else
                        \douzepts\twelveit\fi}
 
   \def\fontetitredeuxb{\ifx\paradouze\oui\onzepts\eleventi\gponze\textfont1=\elevenib\scriptfont1=\nineib
                         \else\douzepts\twelveti\scriptfont1=\twelveib\scriptfont1=\tenib\gpdouze\fi}
                         
\def\fontetitredeuxl{\onzepts\textfont1=\elevenbf\scriptfont1=\ninebf\twelvebf}
  
\def\fontetitretrois{\textfont0=\elevenrm\scriptfont0=\eightrm\textfont1=\eleveni
                      \scriptfont1=\eighti\scriptscriptfont1=\sixi\elevenit}
                      
\def\fontetitrequatre{\textfont0=\elevenrm\scriptfont0=\eightrm\textfont1=\eleveni
                      \scriptfont1=\eighti\scriptscriptfont1=\sixi\elevenrm}
  
  \newcount\titreun\titreun=0
  \newcount\titredeux\titredeux=0
  \newcount\titretrois\titretrois=0
  \newcount\titrequatre\titrequatre=0
  \newcount\enonce\enonce=0
  
  \def\incr#1{\global\advance#1 by 1 {\the #1}}
  \def\avance#1{\global\advance#1 by 1}
  \def\init#1{\global#1=0}
  
  \long\def\Indentation#1#2{\setbox10=\hbox{\fontetitreun#1}
  	                    \ifdim\wd10 < 4mm
                         \setbox10=\hbox to 4mm{\box10\hfill}
                       \else\ifdim\wd10 < 6mm
                         \setbox10=\hbox to 6mm{\box10\hfill}
  	                    \else\ifdim\wd10 < 8mm
                         \setbox10=\hbox to 8mm{\box10\hfill}
                       \else\ifdim\wd10 < 12mm
                         \setbox10=\hbox to 12mm{\box10\hfill}
                       \fi\fi\fi\fi
                       \dimen10=\hsize
                       \advance \dimen10 by -\wd10
                       \noindent \box10 %
                       \ignorespaces
                       \hbox{\vtop{\hsize=\dimen10\raggedright\noindent\fontetitreun#2}}}

  \long\def\paraun#1{\removelastskip\par\medskip\goodbreak\vskip0pt plus.01\vsize\penalty-100
                \vskip0pt plus-.01\vsize
  	              \init{\titredeux}\ifnum\optionparag=1{\init\eqnumber\init\enonce}\else{}\fi
                  \goodbreak{\fontetitreun
  	                \Indentation{\incr{\titreun}.\ }{\fontetitreun #1\par}}\nobreak\medskip}

 %
 %
 \long\def\paraunc#1{\removelastskip\par\bigskip\goodbreak\vskip0pt plus.01\vsize\penalty-100
                \vskip0pt plus-.01\vsize
  	              \init{\titredeux}
                 \ifnum\optionparag=1{\init{\eqnumber}\init\enonce}\else{}\fi
                  \goodbreak
  	                {\parindent0mm\raggedcenter\fontetitreun\incr{\titreun}.\ 
                     \fontetitreun #1\par}\nobreak\medskip}
                     
\newtoks\titreunl
\titreunl={\ifnum\titreun=1{I}\fi%
\ifnum\titreun=2{II}\fi%
\ifnum\titreun=3{III}\fi%
\ifnum\titreun=4{IV}\fi%
\ifnum\titreun=5{V}\fi%
\ifnum\titreun=6{VI}\fi%
\ifnum\titreun=7{VII}\fi%
\ifnum\titreun=8{VIII}\fi%
\ifnum\titreun=9{IX}\fi%
\ifnum\titreun=10{X}\fi%
\ifnum\titreun=11{XI}\fi%
\ifnum\titreun=12{XII}\fi%
\ifnum\titreun=13{XIII}\fi%
}
\long\def\paraunl#1{\removelastskip\par\bigskip\bigskip\goodbreak\vskip0pt plus.01\vsize\penalty-100
                \vskip0pt plus-.01\vsize
  	              \init{\titredeux}\ifnum\optionparag=1{\init\eqnumber\init\enonce}\else{}\fi
                  \goodbreak{\fontetitreunl
  	                \Indentation{\global\advance\titreun by 1{\the\titreunl}.\ }{\fontetitreunl #1\par}}\nobreak\smallskip}

  
  \long\def\paradeux#1{\init{\titretrois}\vskip0pt plus.01\vsize\penalty-10
                \vskip0pt plus-.01\vsize\ifx \elie\oui\medskip\ifnum\titredeux>0\medskip\fi\fi
                 \Indentation{\fontetitredeux\the\titreun${\cdot}$\incr{\titredeux}.}
                              {\fontetitredeux\textfont1=\eleveni#1}\nobreak\par }
  
  \long\def\paradeuxb#1{\init{\titretrois}\vskip0pt plus.001\vsize\penalty-10
                \vskip0pt plus-.01\vsize{\ifx \elie\oui\medskip\ifnum\titredeux>0\medskip\fi\fi
                  \Indentation
  {\fontetitredeuxb\the\titreun${\cdot}$\incr{\titredeux}.}{ \fontetitredeuxb#1}}\nobreak
\smallskip}

\newtoks\titredeuxl
\titredeuxl={\ifnum\titredeux=1{A}\fi%
\ifnum\titredeux=2{B}\fi%
\ifnum\titredeux=3{C}\fi%
\ifnum\titredeux=4{D}\fi%
\ifnum\titredeux=5{E}\fi%
\ifnum\titredeux=6{F}\fi%
\ifnum\titredeux=7{G}\fi%
\ifnum\titredeux=8{H}\fi%
\ifnum\titredeux=9{I}\fi%
\ifnum\titredeux=10{J}\fi%
\ifnum\titredeux=11{K}\fi%
\ifnum\titredeux=12{L}\fi%
\ifnum\titredeux=13{M}\fi%
}
 \long\def\paradeuxl#1{\init{\titretrois}\vskip0pt plus.001\vsize\penalty-10
                \vskip0pt plus-.01
                \vsize \bigskip%
                  \Indentation
     {\fontetitredeuxl\global\advance\titredeux by 1
  \quad \the\titreunl${\cdot}$\the\titredeuxl.}{ \fontetitredeuxl#1}
  \removelastskip\nobreak\smallskip}
  

  \long\def\paratrois#1{\init{\titrequatre}\ifdim\lastskip<\smallskipamount
                \removelastskip\smallskip\fi
                 \vskip0pt plus.01\vsize\penalty-10
                  \vskip0pt
plus-.01\vsize{\ifx \elie\oui\ifnum\titretrois>0\medskip\fi\fi
\Indentation{\fontetitretrois\the\titreun${\cdot}$\the\titredeux${\cdot}$\incr{\titretrois}.\ }
  {\hskip0mm\baselineskip=14pt\fontetitretrois#1}\nobreak\smallskip}}
  
  
  \long\def\paratroisl#1{\init{\titrequatre}\ifdim\lastskip<\smallskipamount
                \removelastskip\fi
                 \vskip0pt plus.01\vsize\penalty-10
                  \vskip0pt
plus-.01\vsize\ifx \elie\oui\bigskip
\fi
\Indentation{\fontetitretrois\quad \quad \the\titreunl{${\cdot}$}\the\titredeuxl${\cdot}$\incr{\titretrois}.\ }
  {\hskip0mm\fontetitretrois#1}\nobreak\smallskip}


  \long\def\paraquatre#1{\ifdim\lastskip<\smallskipamount
                \removelastskip\smallskip\fi
                 \vskip0pt plus.01\vsize\penalty-10
                  \vskip0pt
                  plus-.01\vsize\par
 
\Indentation{\fontetitrequatre \the\titreun{${\cdot}$}\the\titredeux${\cdot}$\the\titretrois${\cdot}$\incr{\titrequatre}.\ }
{\hskip0mm\fontetitrequatre#1}\nobreak\smallskip}


\newtoks\titrequatrel
\titrequatrel={\ifnum\titrequatre=1{a}\fi%
\ifnum\titrequatre=2{b}\fi%
\ifnum\titrequatre=3{c}\fi%
\ifnum\titrequatre=4{d}\fi%
\ifnum\titrequatre=5{e}\fi%
\ifnum\titrequatre=6{f}\fi%
\ifnum\titrequatre=7{g}\fi%
\ifnum\titrequatre=8{h}\fi%
\ifnum\titrequatre=9{i}\fi%
\ifnum\titrequatre=10{j}\fi%
\ifnum\titrequatre=11{k}\fi%
\ifnum\titrequatre=12{l}\fi%
\ifnum\titrequatre=13{m}\fi%
}
\long\def\paraquatrel#1{\ifdim\lastskip<\smallskipamount
                \removelastskip\smallskip\fi
                 \vskip0pt plus.01\vsize\penalty-10
                  \vskip0pt
                  plus-.01\vsize{\bigskip
\Indentation{\global\advance\titrequatre by 1
\fontetitrequatre\quad \quad \quad \the\titreunl${\cdot}$\the\titredeuxl${\cdot}$\the\titretrois${\cdot}$\the\titrequatrel.\ }
{\hskip0mm\fontetitrequatre#1}\nobreak\smallskip}}

\ifx\optionkeys\oui
\def\drefun#1{\definexref{¤#1}{{\the\titreun}}{}} 
\def\drefdeux#1{\definexref{¤#1}{{\the\titreun}.{\the\titredeux}}{}}
\def\dreftrois#1{\definexref{¤#1}{{\the\titreun}.{\the\titredeux}.{\the\titretrois}}{}}
\else
\def\drefun#1{\definexref{prg#1}{{\the\titreun}}{}} 
\def\drefdeux#1{\definexref{prg#1}{{\the\titreun}.{\the\titredeux}}{}}
\def\dreftrois#1{\definexref{prg#1}{{\the\titreun}.{\the\titredeux}.{\the\titretrois}}{}}
\fi

%


  \long\def\propdeux#1#2#3#4{%
       \avance{\enonce}
       \leavevmode\edef\temp{#2}%
         \ifx\temp\empty 
          \else
           \definexref{#2}{#1~{\the\titreun.\the\enonce}}{enonces}
            \definexref{s#2}{{\the\titreun.\the\enonce}}{enonces}
             \fi
\smallskip
      \noindent{\bf#1\ {\bf\the\titreun.\the\enonce{#3}.}\enspace}{\sl#4\par}%
      \ifdim\lastskip<\medskipamount \removelastskip\penalty55\par \fi
   }

  \long\def\propun#1#2#3#4{%
      \avance{\enonce}
       \leavevmode\edef\temp{#2}%
        \ifx\temp\empty 
          \else
           \definexref{#2}{#1~{\the\enonce}}{enonces}
            \definexref{{s#2}}{{\the\enonce}}{enonces}
             \fi
   \par 
     \noindent{\bf#1\ {\bf\the\enonce{#3}.}\enspace}{\sl#4\par}%
     \ifdim\lastskip<\medskipamount \removelastskip\penalty55\medskip\fi
  }
  
  \long\def\prop#1#2#3#4{\ifnum\optionparag=1
                          \propdeux{#1}{#2}{\textfont1=\elevenib#3}{#4} \else\propun{#1}{#2}{\textfont1=\elevenib#3}{#4}\fi}

  \long\def\propt#1#2#3{\ifx\tpf\oui \prop{Th\'eo\-r\`eme}{#1}{#2}{#3}\par
                       \else\prop{Theorem}{#1}{#2}{#3}\par\fi}
  \long\def\Propt#1#2{\propt{#1}{}{#2}}
  \long\def\propl#1#2#3{\ifx\tpf\oui\prop{Lem\-me}{#1}{#2}{#3}\par
                         \else\prop{Lemma}{#1}{#2}{#3}\par\fi}
  
  \long\def\propc#1#2#3{\ifx\tpf\oui\prop{Corol\-laire}{#1}{#2}{#3}\par
                         \else\prop{Corollary}{#1}{#2}{#3}\par\fi}

  \long\def\propd#1#2#3{\ifx\tpf\oui\prop{D\'efi\-nition}{#1}{#2}{#3}\par
                       \else\prop{Definition}{#1}{#2}{#3}\par\fi} 
  
  \long\def\proptd#1#2#3{\ifx\tpf\oui\prop{Th\'eor\`eme et d\'efi\-nition}{#1}{#2}{#3}\par
                       \else\prop{Theorem and definition}{#1}{#2}{#3}\par\fi}


  
  \newcount\optionparag\optionparag=1
  
  \long\def\section#1#2{\ifnum\optionparag=1 \paraun{#2} 
                        \else\goodbreak{\fontetitreun
  	                \Indentation{#1.\ }{#2}}\nobreak\smallskip\fi}

  \def\eqconstruct#1{\ifnum\optionparag=1{\the\titreun\hbox{$\cdot$}#1}\else{#1}\fi}

  
  
  \def\numref{oui}  
  
  \newcount\mesref\mesref=0 
  \def\defbib#1{\ifx\numref\oui\global\advance\mesref by 1 \definexref{#1}{{\the
                 \mesref}}{}\else\definexref{#1}{#1}{}\fi}
  \def\bibtem#1{\defbib{#1}\item{\citer{#1}}}
  \def\citer#1{[\ref{#1}]}

  
  \font\seventeenmsa=msam10 at 17pt    
  \font\fourteenmsa=msam10 at 14pt
  \font\twelvemsa=msam10 at 12pt
  \font\tenmsa=msam10                 
  \font\ninemsa=msam10 at 9pt 
  \font\eightmsa=msam10 at 8pt 
  \font\sevenmsa=msam7 
  \font\sixmsa=msam10 at 6pt
  \font\fivemsa=msam5
  \newfam\msafam\textfont\msafam=\tenmsa\scriptfont\msafam=\sevenmsa\scriptscriptfont\msafam=\fivemsa
  
  \font\seventeenbb=msbm10 at 17pt     
  \font\fourteenbb=msbm10 at 14pt
  \font\twelvebb=msbm10 at 12pt
  \font\tenbb=msbm10                   
  \font\ninebb=msbm10 at 9pt 
  \font\eightbb=msbm10 at 8pt 
  \font\sevenbb=msbm7 
  \font\sixbb=msbm10 at 6pt
  \font\fivebb=msbm5 
  \newfam\bbfam\textfont\bbfam=\tenbb\scriptfont\bbfam=\sevenbb\scriptscriptfont\bbfam=\fivebb
  \def\bb{\fam\bbfam\tenbb}%

  \font\seventeenscaln=eusm10 at 17pt   
  \font\twelvescaln=eusm10 at 12pt
  \font\tenscaln=eusm10                
  \font\ninescaln=eusm10 scaled 900
  \font\eightscaln=eusm10 scaled 800
  \font\sevenscaln=eusm10 scaled 700
  \font\sixscaln=eusm10 scaled 600
   
  \newfam\scalnfam\textfont\scalnfam=\tenscaln\scriptfont\scalnfam=\sevenscaln\scriptscriptfont\scalnfam=\sixscaln
  \def\scaln{\fam\scalnfam\tenscaln}%
  \def\scal{\scaln}
  
  \font\tenscalb=eusb10                

  \font\sevenscalb=eusb10 scaled 700

  \newfam\scalbfam\textfont\scalbfam=\tenscalb\scriptfont\scalbfam=\sevenscalb
  %
  
  %
  %
  \font\fourteenrm=cmr12 scaled 1200
  \font\elevenrm=cmr10 at 11pt
  \font\twelverm=cmr12
  \font\ninerm=cmr9
  \font\eightrm=cmr8      
  \font\sevenrm=cmr7
  \font\sixrm=cmr6

  \font\seventeenpcap=cmcsc10 at 17pt
  \font\tenpcap=cmcsc10                        
  \font\ninepcap=cmcsc9
  \font\eightpcap=cmcsc8
  \font\sevenpcap=cmcsc10 scaled 700
  
  \newfam\pcapfam\textfont\pcapfam=\tenpcap\scriptfont\pcapfam=\sevenpcap
  \def\pcap{\fam\pcapfam\tenpcap}
  
  \font\seventeenrm=cmbx12 scaled 1400

  \font\fourteenbf=cmbx10 scaled 1400
  
  \font\twelvebf=cmbx12
  \font\elevenbf=cmbx10 at 11pt
  \font\ninebf=cmbx9  
  \font\eightbf=cmbx8
  \font\sixbf=cmbx6
  
  \font\tengot=eufm10                           
   
  \font\eightgot=eufm10 at 8truept 
  \font\sevengot=eufm7 
  \font\sixgot=eufm10 at 6 truept 
   
  \newfam\gotfam
  \textfont\gotfam=\tengot\scriptfont\gotfam=\sevengot\scriptscriptfont\gotfam=\sixgot
  %

  
  \def\tit{%
  \textfont0=\seventeenrm\scriptfont0=\tenrm\def\rm{\fam0\seventeenrm}%
  \textfont1=\seventeenib\scriptfont1=\twelveib%
  \textfont2=\seventeensy\scriptfont2=\twelvesy\scriptscriptfont2=\ninesy
  \textfont3=\seventeenex
  \textfont\itfam=\seventeenti
  \def\it{\fam\itfam\seventeenti}%
  \textfont\bbfam=\seventeenbb \scriptfont\bbfam=\twelvebb
  \def\bb{\fam\bbfam\seventeenbb}%
  \textfont\msafam=\seventeenmsa\scriptfont\msafam=\twelvemsa
  \textfont\scalnfam=\seventeenscaln
  \def\pcap{\fam\pcapfam\seventeenpcap}
  \normalbaselineskip=25pt\normalbaselines\rm}

  \font\seventeenti=cmbxti10 scaled 1680
  
  \font\fourteenti=cmbxti10 at 14pt
  
  \font\twelveti=cmbxti10 scaled 1200
  \font\eleventi=cmbxti10 at 11pt

  %
  %
  \font\twelveit=cmti12	
  \font\elevenit=cmti10 scaled 1100
  \font\nineit=cmti9
  \font\eightit=cmti8
  \font\sevenit=cmti7

  %
  %
 
 \font\seventeenib=cmmib10 scaled 1680
  \font\fourteenib=cmmib10 scaled 1400
  \font\twelveib=cmmib10 scaled 1200
  \font\elevenib=cmmib10 scaled 1100
  \font\tenib=cmmib10
\font\eightib=cmmib10 scaled 800
  \font\nineib=cmmib10 scaled 900
\font\sevenib=cmmib10 scaled 700
\font\sixib=cmmib10 scaled 600
\font\fiveib=cmmib10 scaled 500

\ifx\ITAN\oui
\else
\innernewfam\cmmibfam
\textfont\cmmibfam=\tenib
\scriptfont\cmmibfam=\sevenib
\scriptscriptfont\cmmibfam=\fiveib
\def\ib{\fam\cmmibfam\tenib}
\fi

  %
  %
  
  \font\eleveni=cmmi10 scaled 1100
  \font\ninei=cmmi9
  \font\eighti=cmmi8 
  \font\seveni=cmmi7 			                
  \font\sixi=cmmi6
  
  \font\ninesl=cmsl9                    
  \font\eightsl=cmsl8 
  \font\sevensl=cmsl10 at 7pt

  \font\ninett=cmtt9                    
  \font\eighttt=cmtt8
  \font\seventt=cmtt10 scaled 700

  \font\seventeensy=cmsy10 scaled 1680    
  \font\fourteensy=cmsy10 scaled 1400
  \font\twelvesy=cmsy10 scaled 1176
  
  \font\ninesy=cmsy9                      
  \font\eightsy=cmsy8
  \font\sixsy=cmsy6
  \font\seventeenex=cmex10 at 17pt
  \font\fourteenex=cmex10 at 14pt
  \font\twelveex=cmex10 at 12pt
  \font\nineex=cmex10 at 9pt
  \font\eightex=cmex10 at 8pt
  \font\sevenex=cmex10 at 7pt
  \font\sixex=cmex10 at 6pt
  \font\fiveex=cmex10 at 5pt
  
   
  \font\fourteengp=cmmi10 at 14pt
  
  \font\twelvegp=cmmib10 at 12pt
  \font\elevengp=cmmib10 at 11pt
  \font\tengp=cmmib10                          
  \font\ninegp=cmmib10 at 9pt 
  \font\eightgp=cmmib8 
   
  \font\sixgp=cmmib6


  \def\gponze{\textfont0=\elevenbf\scriptfont0=\eightbf\scriptscriptfont0=\sixbf
           \textfont1=\elevengp\scriptfont1=\eightgp\scriptscriptfont1=\sixgp}
  \def\gpdouze{\textfont0=\twelvebf\scriptfont0=\tenbf\scriptscriptfont0=\ninebf
           \textfont1=\twelvegp\scriptfont1=\tengp\scriptscriptfont1=\ninegp}        
  
 \def\gpquatorze{\textfont0=\fourteenbf\scriptfont0=\twelvebf\scriptscriptfont0=\elevenbf
           \textfont1=\fourteengp\scriptfont1=\twelvegp\scriptscriptfont1=\elevengp}

  
  \expandafter\chardef\csname pre amssym.def at\endcsname=\the\catcode`\@
  \catcode`\@=11
  \def\undefine#1{\let#1\undefined}
  \def\newsymbol#1#2#3#4#5{\let\next@\relax
   \ifnum#2=\@ne\let\next@\msafam@\else
   \ifnum#2=\tw@\let\next@\bbfam@\fi\fi
   \mathchardef#1="#3\next@#4#5}
  \def\mathhexbox@#1#2#3{\relax
   \ifmmode\mathpalette{}{\m@th\mathchar"#1#2#3}%
   \else\leavevmode\hbox{$\m@th\mathchar"#1#2#3$}\fi}
  \def\hexnumber@#1{\ifcase#1 0\or 1\or 2\or 3\or 4\or 5\or 6\or 7\or 8\or
   9\or A\or B\or C\or D\or E\or F\fi}
  
  \def\setboxz@h{\setbox\z@\hbox}
  \def\wdz@{\wd\z@}
  \def\boxz@{\box\z@}
  
  \edef\msafam@{\hexnumber@\msafam}
  \mathchardef\dabar@"0\msafam@39
  
  \edef\bbfam@{\hexnumber@\bbfam}
  \def\widehat#1{\setboxz@h{$\m@th#1$}%
   \ifdim\wdz@>\tw@ em\mathaccent"0\bbfam@5B{#1}%
   \else\mathaccent"0362{#1}\fi}
  \def\widetilde#1{\setboxz@h{$\m@th#1$}%
   \ifdim\wdz@>\tw@ em\mathaccent"0\bbfam@5D{#1}%
   \else\mathaccent"0365{#1}\fi}
  \newsymbol\leqq 1335          
  \newsymbol\leqslant 1336
  \newsymbol\lessgtr 1337       
  \newsymbol\backprime 1038     
  \newsymbol\risingdotseq 133A  
  \newsymbol\fallingdotseq 133B 
  \newsymbol\succcurlyeq 133C   
  \newsymbol\geqq 133D          
  \newsymbol\geqslant 133E
  \newsymbol\nmid 232D
  \newsymbol\nexists 2040
  \newsymbol\smallsetminus 2272
  \newsymbol\varnothing 203F
  
  \catcode`\@=\active

  \catcode`\@=11
  \newcount\typofr\typofr=1
  
  \catcode`\;=\active
  \def;{\ifnum\typofr=1\relax\ifhmode\ifdim\lastskip>\z@\unskip\fi
     \kern.2em\fi\string;\else\string;\fi}
  
  \catcode`\:=\active
  \def:{\ifnum\typofr=1\relax\ifhmode\ifdim\lastskip>\z@\unskip\fi
  \penalty\@M\ \fi\string:\else\string:\fi}
  
  \catcode`\!=\active
  \def!{\ifnum\typofr=1\relax\ifhmode\ifdim\lastskip>\z@\unskip\fi
     \kern.2em\fi\string!\else\string!\fi}
  
  \catcode`\?=\active
  \def?{\ifnum\typofr=1\relax\ifhmode\ifdim\lastskip>\z@\unskip\fi
     \kern.2em\fi\string?\else\string?\fi}

  \def\francais{\typofr=1\def\tpf{oui}}
  
  \def\oui{oui}
  \francais
  
  \catcode`\@=12
  

\ifx\textures\oui
\def\raye #1|{\leavevmode\setbox1=\hbox{#1}%
\raise .5pt\hbox to \wd1{\xleaders\hbox{\rge{$ \sct / $}%
\kern 1pt}\hfill\kern -1pt }\kern -\wd1 \unhbox1\relax }

\def\barre #1|{\leavevmode\setbox1=\hbox{#1}%
\rlap{\Red\vrule height 2.4pt depth -1.2pt width \wd1}\Black \unhbox1\relax}
\else
\def\raye #1|{\leavevmode\setbox1=\hbox{#1}%
\raise .5pt\hbox to \wd1{\xleaders\hbox{\rge{$ \sct / $}%
\kern 1pt}\hfill\kern -1pt }\kern -\wd1 \unhbox1\relax }

\def\barre #1|{\leavevmode\setbox1=\hbox{#1}%
\rlap{\color{red}\vrule height 2.4pt depth -1.2pt width \wd1}\color{black} \unhbox1\relax}

\fi
  

  
  \def\og{\leavevmode\raise.24ex\hbox{$\scriptscriptstyle\langle\!\langle\>$}}    
  \def\fg{\leavevmode\raise.24ex\hbox{$\scriptscriptstyle\>\rangle\!\rangle$}}    
  \def\chv#1{\left\langle#1\right\rangle}

  \def\dd{{\rm d}}

  \def\z{{\bb Z}}
  \def\r{{\bb R}}

  \def\PP{{\bb P}}

  \def\B{{\scal B}}

  \def\E{{\scal E}}

  \def\M{{\scal M}}
  
  \def\O{{\scal O}}
  \def\P{{\scaln P}}

  \def\frac#1#2{{#1\over #2}}
  \def\di#1#2{\sct#1\atop{\sct#2}}

  \def\theoreme#1{\rm th.\thinspace#1}
  

  \def\¤{\S\thinspace}

  \def\¥{$\bullet$ }
  
  
  \def\e{{\rm e}}

  \def\no#1{\Vert#1\Vert}
  
  \def\NO#1{\left\Vert#1\right\Vert}

  \def\epsilon{\varepsilon}

  \def\phi{\varphi}
  \def\theta{\vartheta}
  \def\rho{\varrho}
  \def\dm{{\textstyle{1\over 2}}}
  \def\txt{\textstyle}
  \def\dsp{\displaystyle}
  \def\sct{\scriptstyle}
  \def\pf{\noi{\it Proof. }}
  \def\nid{\ifnum\typofr=1\par\noindent{\it D\'emonstration. }\else\pf\fi}
  \def\noi{\noindent}
  \def\rem{\ifnum\typofr=1\noi{\it Remarque.}\ \else\noi{\it Remark.}\ \fi}
  \def\rems{\ifnum\typofr=1\noi{\it Remarques.}\ \else\noi{\it Remarks.}\ \fi}

  \def\sset{\smallsetminus}

  \def\1{{\bf 1}}
  \def\|{\Vert}

  \def\leq{\leqslant}
  \def\geq{\geqslant}
  
  \def\cf{{cf.}}
  
  \def\eg{{e.g.}}
  

  \def\fl#1{\left\lfloor #1 \right\rfloor}

  \def\log{\mathop{\rm log}\nolimits}
  \def\ft#1#2{{\txt{#1\over #2}}}



\def\Vbs#1{\bigg|#1\bigg|}


  \def\pmb#1{\setbox0=\hbox{#1}%
  \kern-.025em\copy0\kern-\wd0\kern.05em\copy0\kern-\wd0\kern-.025em\raise .0433em\box0 }

  
  \skewchar\eighti='177 \skewchar\sixi='177
  \skewchar\eightsy='60 \skewchar\sixsy='60
  
  \def\eightpoint{%
  \textfont0=\eightrm\scriptfont0=\sixrm\scriptscriptfont0=\fiverm
  \def\rm{\fam0\eightrm}%
  \textfont1=\eighti\scriptfont1=\sixi
  \scriptscriptfont1=\fivei\def\oldstyle{\fam1\seveni}%
  \textfont2=\eightsy\scriptfont2=\sixsy\scriptscriptfont2=\fivesy
  \textfont3=\eightex\scriptfont3=\sixex
  \textfont\itfam=\eightit
  \def\it{\fam\itfam\eightit}%
  \textfont\slfam=\eightsl
  \def\sl{\fam\slfam\eightsl}%
  \textfont\bbfam=\eightbb \scriptfont\bbfam=\sixbb\scriptscriptfont\bbfam=\fivebb
  \def\bb{\fam\bbfam\eightbb}%
  \textfont\msafam=\eightmsa\scriptfont\msafam=\sixmsa
  \textfont\scalnfam=\eightscaln
  \def\scaln{\fam\scalnfam\eightscaln}
  \textfont\ttfam=\eighttt
  \def\tt{\fam\ttfam\eighttt}%
\textfont\gotfam=\eightgot
  \textfont\bffam=\eightbf\scriptfont\bffam=\sixbf\scriptscriptfont\bffam=\fivebf
  \def\bf{\fam\bffam\eightbf}%
  \ifx\ITAN\oui\else\textfont\cmmibfam=\eightib
       \scriptfont\cmmibfam=\sixib
        \scriptscriptfont\cmmibfam=\fiveib
         \def\ib{\fam\cmmibfam\eightib}
   \fi
  \textfont\pcapfam=\eightpcap
  \def\pcap{\fam\pcapfam\eightpcap}
  \abovedisplayskip=2pt plus2pt minus 2pt
  \belowdisplayskip=2pt plus1pt minus 2pt
  \abovedisplayshortskip= 1pt plus 2pt minus 1pt
  \belowdisplayshortskip= 1pt plus 2pt minus 1pt
  \smallskipamount=2pt plus 1pt minus 2pt
  \medskipamount=3pt plus 2pt minus 2pt
  \bigskipamount=7pt plus 3pt minus 3pt
  \setbox\strutbox=\hbox{\vrule height 5pt depth 2pt width 0pt}%
  \normalbaselineskip=9pt\normalbaselines\rm}

  \def\({\left(}
  \def\){\right)}
  
  \def\footnoterule{\kern -2pt\hrule width 7truecm\kern 2.4pt}
  
  \def\xnotedef#1{\definexref{#1}{\noexpand\number\footnotenumber}{Note}}%

  
  
  \def\ninepoint{%
  \textfont0=\ninerm\scriptfont0=\sixrm\scriptscriptfont0=\fiverm
  \def\rm{\fam0\ninerm}%
  \textfont1=\ninei\scriptfont1=\sixi
  \scriptscriptfont1=\fivei\def\oldstyle{\fam1\ninei}%
  \textfont2=\ninesy\scriptfont2=\sixsy\scriptscriptfont2=\fivesy
  \textfont3=\nineex\scriptfont3=\sixex
  \textfont\itfam=\nineit
  \def\it{\fam\itfam\nineit}%
  \textfont\slfam=\ninesl
  \def\sl{\fam\slfam\ninesl}%
  \textfont\bbfam=\ninebb\scriptfont\bbfam=\sixbb\scriptscriptfont\bbfam=\fivebb
  \def\bb{\fam\bbfam\ninebb}%
  \textfont\msafam=\ninemsa\scriptfont\msafam=\sixmsa\scriptscriptfont\msafam=\fivemsa
  \textfont\scalnfam=\ninescaln
  \def\scaln{\fam\scalnfam\ninescaln}
  \textfont\ttfam=\ninett
  \def\tt{\fam\ttfam\ninett}%
  \textfont\bffam=\ninebf\scriptfont\bffam=\sixbf\scriptscriptfont\bffam=\fivebf
  \def\bf{\fam\bffam\ninebf}%
  \abovedisplayskip=3pt plus2pt minus 2pt
  \belowdisplayskip=3pt plus1pt minus 2pt
  \abovedisplayshortskip= 2pt plus 2pt minus 1pt
  \belowdisplayshortskip= 2pt plus 2pt minus 1pt
  \smallskipamount=2pt plus 1pt minus 2pt
  \medskipamount=3pt plus 2pt minus 2pt
  \bigskipamount=7pt plus 3pt minus 3pt
  \setbox\strutbox=\hbox{\vrule height 5pt depth 2pt width 0pt}%
  \normalbaselineskip=10.5pt plus.3pt minus.3pt\normalbaselines\rm}

  \def\sevenpoint{%
  \textfont0=\sevenrm\scriptfont0=\sixrm\scriptscriptfont0=\fiverm
  \def\rm{\fam0\sevenrm}%
  \textfont1=\seveni\scriptfont1=\sixi
  \scriptscriptfont1=\fivei\def\oldstyle{\fam1\seveni}%
  \textfont2=\sevensy\scriptfont2=\sixsy\scriptscriptfont2=\fivesy
  \textfont3=\sevenex\scriptfont3=\fiveex
  \textfont\itfam=\sevenit
  \def\it{\fam\itfam\sevenit}%
  \textfont\slfam=\sevensl
  \def\sl{\fam\slfam\sevensl}%
  \textfont\bbfam=\sevenbb \scriptfont\bbfam=\sixbb\scriptscriptfont\bbfam=\fivebb
  \def\bb{\fam\bbfam\sevenbb}%
  \textfont\msafam=\sevenmsa\scriptfont\msafam=\sixmsa
  \textfont\scalnfam=\sevenscaln
  \def\scaln{\fam\scalnfam\sevenscaln}
  \textfont\bffam=\sevenbf\scriptfont\bffam=\sixbf\scriptscriptfont\bffam=\fivebf
  \def\bf{\fam\bffam\sevenbf}%
  \textfont\ttfam=\seventt
  \abovedisplayskip=2pt plus2pt minus 2pt
  \belowdisplayskip=2pt plus1pt minus 2pt
  \abovedisplayshortskip= 1pt plus 2pt minus 1pt
  \belowdisplayshortskip= 1pt plus 2pt minus 1pt
  \smallskipamount=2pt plus 1pt minus 2pt
  \medskipamount=3pt plus 2pt minus 2pt
  \bigskipamount=7pt plus 3pt minus 3pt
  \setbox\strutbox=\hbox{\vrule height 5pt depth 2pt width 0pt}%
  \normalbaselineskip=9pt\normalbaselines\rm}

 \def\onzepts{%
 \textfont0=\elevenrm\scriptfont0=\ninerm
 \textfont1=\elevenib\scriptfont1=\ninei}

\def\douzepts{%
  \textfont0=\twelverm\scriptfont0=\tenrm\def\rm{\fam0\twelverm}%
  \textfont1=\twelveib\scriptfont1=\teni%
  \textfont2=\twelvesy\scriptfont2=\tensy\scriptscriptfont2=\eightsy
  \textfont3=\twelveex
  \textfont\itfam=\twelveti
  \def\it{\fam\itfam\twelveti}%
  \textfont\bffam=\twelvebf\scriptfont\bffam=\tenbf\scriptscriptfont\bffam=\eightbf
  \def\bf{\fam\bffam\twelvebf}%
  \textfont\bbfam=\twelvebb \scriptfont\bbfam=\tenbb
  \def\bb{\fam\bbfam\twelvebb}%
  \textfont\msafam=\twelvemsa\scriptfont\msafam=\tenmsa
  \textfont\scalnfam=\twelvescaln
  \normalbaselineskip=15pt\normalbaselines\rm}

\def\quatorzepts{%
  \textfont0=\fourteenrm\scriptfont0=\twelverm\def\rm{\fam0\fourteenrm}%
  \textfont1=\fourteenib\scriptfont1=\twelveib%
  \textfont2=\fourteensy\scriptfont2=\twelvesy\scriptscriptfont2=\tensy
  \textfont3=\fourteenex
  \textfont\itfam=\fourteenti
  \def\it{\fam\itfam\fourteenti}%
  \textfont\bffam=\fourteenbf\scriptfont\bffam=\twelvebf\scriptscriptfont\bffam=\tenbf
  \def\bf{\fam\bffam\fourteenbf}%
  \textfont\bbfam=\fourteenbb \scriptfont\bbfam=\twelvebb
  \def\bb{\fam\bbfam\fourteenbb}%
  \textfont\msafam=\fourteenmsa\scriptfont\msafam=\twelvemsa
  \textfont\scalnfam=\twelvescaln
  \normalbaselineskip=18pt\normalbaselines\rm}


\def\AA{{\it Acta Arith.}}

\def\picture #1 by #2 (#3){\leavevmode\vbox to #2{
     \hrule width #1 height 0pt depth 0pt
      \vfill
       \special{picture #3}}}

\def\illustration #1 by #2 (#3) scaled #4{\dimen1=#2
  \divide\dimen1 by 1000
  \multiply\dimen1 by #4
  \vtop to \dimen1{\dimen1=#1
  \divide\dimen1 by 1000
  \multiply\dimen1 by #4
  \hsize=\dimen1\vss
  \noindent\special{illustration #3 scaled #4}}}

\beginpackages
\usepackage{color}
\endpackages 
\long\def\rge#1{{\color{red}#1}}

 \fi

\dimstand  
 \francais

\ifx\optionkeymacros\undefined\else \fi

\catcode`\Œ=\active\defŒ{{\aa}}       
\catcode`\º=\active\defº{\int}        
\catcode`\=\active\def{\c c}        
\catcode`\¶=\active\def¶{\partial}    
\catcode`\Ä=\active\defÄ{\oint}       
\catcode`\Æ=\active\defÆ{\triangle}   
\catcode`\Â=\active\defÂ{\neg}        
\catcode`\µ=\active\defµ{\mu}         
\catcode`\¿=\active\def¿{{\o}}        
\catcode`\¹=\active\def¹{\pi}         
\catcode`\Ï=\active\defÏ{{\oe}}       
\catcode`\§=\active\def§{{\ss}}       
\catcode`\ =\active\def {\dagger}     
\catcode`\Ã=\active\defÃ{\sqrt}       
\catcode`\·=\active\def·{\Sigma}      
\catcode`\Å=\active\defÅ{\approx}     
\catcode`\½=\active\def½{\Omega}      
\catcode`\£=\active\def£{{\it\$}}     
\catcode`\°=\active\def°{\infty}      
\catcode`\¤=\active\def¤{{\S}}        
\catcode`\¦=\active\def¦{{\P}}        
\catcode`\¥=\active\def¥{\bullet}     
\catcode`\»=\active\def»{\leavevmode\raise.585ex\hbox{\b a}}      
\catcode`\¼=\active\def¼{\leavevmode\raise.6ex\hbox{\b o}}        
\catcode`\­=\active\def­{\not=}       
\catcode`\²=\active\def²{\leq}        
\catcode`\³=\active\def³{\geq}        
\catcode`\Ö=\active\defÖ{\div}        
\catcode`\É=\active\defÉ{{\dots}}     
\catcode`\¾=\active\def¾{{\ae}}       
\catcode`\Ç=\active\defÇ{\og}         
\catcode`\Ò=\active\defÒ{``}          
\catcode`\Á=\active\defÁ{!`}          
\catcode`\¢=\active\def¢{\rlap/c}     
\catcode`\Ô=\active\defÔ{`}           
\catcode`\Õ=\active\defÕ{'}           


\catcode`\=\active\def{{\AA}}       
\catcode`\'=\active\def'{\c C}        
\catcode`\¯=\active\def¯{{\O}}        
\catcode`\¸=\active\def¸{\Pi}         
\catcode`\Î=\active\defÎ{{\OE}}       
\catcode`\®=\active\def®{{\AE}}       
\catcode`\×=\active\def×{\diamond}    
\catcode`\¡=\active\def¡{\accent'27}  
\catcode`\Ó=\active\defÓ{''}          
\catcode`\±=\active\def±{\pm}         
\catcode`\È=\active\defÈ{\fg}         
\catcode`\À=\active\defÀ{?`}          
\catcode`\Ð=\active\defÐ{--}          
\catcode`\Ñ=\active\defÑ{---}         


\catcode`\Š=\active\defŠ{\"a}        
\catcode`\'=\active\def'{\"e}        
\catcode`\•=\active\def•{\"{\i}}     
\catcode`\š=\active\defš{\"o}        
\catcode`\Ÿ=\active\defŸ{\"u}        
\catcode`\Ø=\active\defØ{\"y}        
\catcode`\å=\active\defå{\^A}        
\catcode`\€=\active\def€{\"A}        
\catcode`\…=\active\def…{\"O}        
\catcode`\†=\active\def†{\"U}        
\catcode`\‡=\active\def‡{\'a}        
\catcode`\Ž=\active\defŽ{\'e}        
\catcode`\'=\active\def'{\'{\i}}     
\catcode`\—=\active\def—{\'o}        
\catcode`\œ=\active\defœ{\'u}        
\catcode`\ƒ=\active\defƒ{\'E}        
\catcode`\æ=\active\defæ{\^E}        
\catcode`\é=\active\defé{\`E}        %
\catcode`\ˆ=\active\defˆ{\`a}        
\catcode`\=\active\def{\`e}        
\catcode`\"=\active\def"{\`{\i}}     
\catcode`\˜=\active\def˜{\`o}        
\catcode`\=\active\def{\`u}        
\catcode`\Ë=\active\defË{\`A}        
\catcode`\‹=\active\def‹{\~a}        
\catcode`\–=\active\def–{\~n}        
\catcode`\›=\active\def›{\~o}        
\catcode`\Ì=\active\defÌ{\~A}        
\catcode`\"=\active\def"{\~N}        
\catcode`\Í=\active\defÍ{\~O}        
\catcode`\‰=\active\def‰{\^a}        
\catcode`\=\active\def{\^e}        
\catcode`\"=\active\def"{\^{\i}}     
\catcode`\™=\active\def™{\^o}        
\catcode`\ž=\active\defž{\^u}        

\let\optionkeymacros\null

\def\paradouze{oui}


\def\voircorrections{oui} 

\ifx\montrerlabels\oui\input montrerlabels.tex\fi

%
%

\newcount\paras\paras=0
\newcount\sparas\sparas=0
\def\tocsectionentry#1#2{\init{\sparas}\avance\paras
                          {\quad\bf\the\paras\quad }{\hskip-2mm
#1\dotfill\hskip3mm\rm#2}\par}%
\def\tocsubsectionentry#1#2{\avance\sparas
                          {\qquad\eightpoint\it\the\paras.\the\sparas}
{\hskip-2mm\eightpoint\it#1\dotfill\hskip3mm\rm#2}\par}%

%
%

\def\??{\leftnote{{??}}}

\ifx\voircorrections\oui\else\def\rge#1{#1}\fi 

\dateheure


\font\tenib=cmmib10
            \font\nineib=cmmib10 scaled 900
\font\sevenib=cmmib10 scaled 700
\font\eightib=cmmib10 scaled 800
\font\fiveib=cmmib10 scaled 500

\def\ITAN{oui}

    \font\tenrsfs=rsfs7 at 10pt

    \font\sevenrsfs=rsfs7
    \font\sixrsfs=rsfs7 at 6pt
    
\newfam\rsfsfam\textfont\rsfsfam=\tenrsfs\scriptfont\rsfsfam=\sevenrsfs\scriptscriptfont\rsfsfam=\sixrsfs

\def\VV{{\bb V}}
\def\EE{{\bb E}}

\def\fl#1{\left\lfloor #1 \right\rfloor}

\def\titrart{Somme des chiffres et changement de base}
\def\auteurs{RŽgis de la Bret\`eche, Thomas Stoll \& GŽrald Tenenbaum}
\hautspages{\auteurs}{\titrart}
\titrecentre{\titrart}
\bigskip
\centerline{\auteurs}
\bigskip\bigskip
{\eightpoint\leftskip1cm\rightskip1cm
\noi{\bf Abstract.} 
 Let $s_a(n)$ denote the sum of digits of an integer $n$ in the base $a$ expansion. Answering, in an extended form, a question of Deshouillers, Habsieger,  Laishram, and Landreau, we show that, provided $a$ and $b$ are multiplicatively independent, any positive real number is a limit point of the sequence $\{s_b(n)/s_a(n)\}_{n=1}^{\infty}$. We also provide upper and lower bounds for the counting functions of the corresponding subsequences. \par 
\PAR
\medskip\noi
{\bf Keywords:} sum of digits, multiplicative independence, exponent of irrationality, binomial recentering.
\medskip\noi
{\bf AMS 2010 Classification:} 11A63, 11K16, 11K60, 11J82.\par}
\par 
\bigskip

\medskip 
\paraun{Introduction et ŽnoncŽ des rŽsultats}
RŽpondant ˆ une question de Steinhaus, Cassels \citer{Ca59} a construit en 1959 un ensemble de mesure de Hausdorff positive dont les ŽlŽments sont normaux en toute base qui n'est pas une puissance de $3$. Dans le mme esprit, alors qu'il est gŽnŽralement conjecturŽ que les dŽveloppements des entiers dans des bases multiplicativement indŽpendantes sur $\z$ sont statistiquement indŽpendants, il est naturel d'attendre un phŽnomne de dŽpendance pour une infinitŽ d'entiers exceptionnels. Deshouillers {\it et al.} \citer{DHLL17}  ont ainsi observŽ numŽriquement une co•ncidence des sommes des chiffres $s_2(n)$ et $s_3(n)$ en bases $2$ et $3$ pour une proportion significative d'entiers $n$ et Žtabli que $|s_3(n)-s_2(n)|\leqslant \delta\log n$ pour une infinitŽ d'entiers $n$ avec $\delta\approx0,14572$. Ils posent alors le problme de l'existence d'une suite infinie d'entiers $n$ pour lesquels $|s_2(n)-s_3(n)|$ est ``significativement petit''. Au vu du rŽsultat prŽcŽdemment citŽ, il est naturel de se demander si l'on a $s_2(n)\sim s_3(n)$ pour une sous-suite convenable des entiers naturels. 
Nous nous proposons ici de rŽpondre positivement ˆ cette question.
\par
Notre approche fonctionne pour les sommes des chiffres $s_a(n)$ et $s_b(n)$ en bases respectives $a$ et $b$ lorsque $\log a$ et $\log b$ sont linŽairement indŽpendants sur $\z$. La preuve nŽcessite seulement que $\vartheta:=(\log a)/\log b$ possde un exposant d'irrationalitŽ fini, ce qui est toujours le cas, d'aprs~\citer{B72}, si le rapport est irrationnel: cela signifie qu'il existe $\gamma\geqslant 2$ tel que
$$\Big|\vartheta-{r\over q}\Big|\gg {1\over q^{\gamma}}\qquad \big(q\geqslant 1,\,(r,q)=1\big).\eqdef{defdelta}$$  
AmŽliorant un rŽsultat de \citer{S07}, Wu et Wang \citer{WW14} ont montrŽ que $\gamma=5,117$ est admissible pour $\vartheta:=(\log 2)/\log 3$. 
\par 
Nous Žtablissons en fait que, pour tout $\tau>0 $, il existe une infinitŽ d'entiers~$n$ tels que
$$s_b(n)\sim \tau s_a(n).\eqdef{sbtsa}$$
\par 
Posons $$ \tau_0=\tau_0(a,b):={(b-1)\log a\over (a-1)\log b}\cdot\eqdef{deftab}$$
Lorsque $\tau=\tau_0$  la relation \eqref{sbtsa} a trivialement lieu sur une suite de densitŽ unitŽ. En effet,  pour presque tout entier $n$, nous avons
$$s_b(n)\sim {b-1\over  2\log b}\log n,\qquad  s_a(n)\sim {a-1\over  2\log a}\log n.$$
\par  
Pour $\tau>0$, dŽfinissons $\varrho_1=\varrho_1(\tau):=\tau/(2\tau_0)$, $\varrho_2=\varrho_2(\tau):=\tau_0/(2\tau)$,
$$\normalbaselineskip=30pt\cases{\dsp c_1=c_1(a,b;\tau):= {\varrho_1\over \log b}\log \Big({1\over \varrho_1}\Big)+{1-\varrho_1\over \log b}\log \Big({1\over 1-\varrho_1}\Big)  & si $\tau<  2\tau_0$,\cr\dsp c_2=c_2(a,b;\tau):=
 {\varrho_2\over \log a}\log \Big({1\over \varrho_2}\Big)+{1-\varrho_2\over \log a}\log \Big({1\over 1-\varrho_2}\Big)  & si  $\tau>\dm\tau_0,$\cr} $$ et
$$ c_0=c_0(a,b;\tau):=\normalbaselineskip=13pt\cases{c_1 & si $\tau\leqslant \dm\tau_0$,\cr  \dsp \max(c_1, c_2) & si $\dm \tau_0<\tau<2\tau_0$,
\cr  \dsp c_2 & si $\tau\geqslant 2\tau_0$.\cr} $$
\goodbreak\par \goodbreak
Ainsi, lorsque $\tau=1$, $a=2$, $b=3$, nous avons $$\tau_0=(\log 4)/\log 3\approx 1,26186, \quad\varrho_2=(\log 2)/\log 3, \quad c_0= c_2\approx 0,94996. $$
\par 
Nous obtenons le rŽsultat suivant,  dans lequel nous notons $M:=\max(a,b)$ et 
$$\eqalign {\Lambda:=\cases{\varrho_1(1-\varrho_1) & si $\tau\leqslant\tau_0$\cr
\varrho_2(1-\varrho_2)& si $\tau> \tau_0$}.\cr}$$

\Propt{ths2s3}{Soient $\tau>0 $, $a$, $b$ deux entiers $\geqslant 2$ multiplicativement indŽpendants, et  $c<c_0(a,b;\tau)$. Lorsque $x$ tend vers l'infini, il existe $\gg x^c$ entiers positifs~$n$ n'excŽdant pas~$x$ tels que $s_b(n)\sim \tau s_a(n)$  et, plus prŽcisŽment,   si $\gamma$ est un exposant d'irrationalitŽ de $\vartheta:=(\log  a)/\log  b$,     
$$s_b(n)=\tau s_a(n)\bigg\{1+O\bigg({1\over (\log n)^{\sigma/\gamma}}\bigg)\bigg\},\eqdef{majsa-tsb}$$
 pour tout $\sigma\in]0, \Lambda/(6M^3\log M) [$.\par 
De plus, pour tout nombre rŽel $\tau\neq\tau_0(a,b)$, il existe un exposant $d_0(\tau)=d_0(\tau;a,b)<1$ tel que, lorsque $x\to\infty$, la relation asymptotique \eqref{sbtsa} soit rŽalisŽe pour au plus $\ll x^{d_0(\tau)+o(1)}$ entiers positifs~$n$ n'excŽdant pas $x$.} 
\medskip
\par 
Nous explicitons au paragraphe \ref{prgmaj} une valeur admissible de $d_0(\tau;a,b)$ pour chaque \hbox{$\tau\neq\tau_0$}. Ainsi, $d_0(1;2,3)\approx 0.993702$. Ë titre de comparaison, mentionnons qu'une hypothse d'indŽpen\-dance des reprŽsentations en bases 2 et 3 fournit un cardinal $x^{t+o(1)}$ avec $t\approx 0,981513$.
\par 
 Nous n'avons pas cherchŽ ˆ optimiser la valeur de $\sigma$ dans l'ŽnoncŽ du \ref{ths2s3}. 
\par 
Il est par ailleurs ˆ noter que, dans le cas $\tau=1$,  le problme de l'existence d'une suite infinie satisfaisant \eqref{sbtsa} est trivial lorsque $\log a$ et $\log b$ sont linŽairement dŽpendants, disons $a^r=b^s$ avec $r$ et $s$ entiers positifs. En effet, il suffit alors de considŽrer les entiers $n$ dont le dŽveloppement dans la base $a^r$ ne comporte que des chiffres 0 ou 1, pour lesquels nous avons $s_a(n)=s_b(n)$. Nous pouvons donc finalement Žnoncer que l'Žquivalence asymptotique $s_a(n)\sim s_b(n)$ a lieu pour tous $a,\,b\geqslant 2$ pour une suite de fonction de comptage $\gg x^c$ o $c=c(a,b)>0$.
\par 
Terminons cette introduction par une remarque mŽthodologique relative aux liens de ce problme avec la thŽorie des suites automatiques. Pour toute base $a$, les ensembles de niveau $E_a(k):=\{n\geqslant 1:s_a(n)=k\}$ sont en effet $a$-automatiques \citer{AS03}, ou encore $a$-reconnaissables, au sens o l'ensemble des mots sur l'alphabet $\{0, 1,\ldots , a-1\}$  repr\'esentant les \'el\'ements de $E_a(k)$ est, pour chaque $k$, une partie reconnaissable par automate du mono•de libre $\{0, 1, \ldots, a-1\}^*$. \par 
 Or, le thŽorme de Cobham \citer{Cob69} fournit un critre pour qu'un ensemble $E$ d'entiers soit simultanŽment reconnaissable  en bases $a$ et $b$ lorsque $a$ et $b$ sont multiplicativement indŽpendants: il en est ainsi, si, et seulement si, $E$ est reprŽsentable comme l'union d'un ensemble fini et d'un nombre fini de progressions arithmŽtiques. Ainsi, sous l'hypothse indiquŽe, $E_a(k)$ n'est pas $b$-reconnaissable. Il s'ensuit que les rŽsultats de rŽpartition connus pour les suites automatiques (\eg, relatifs ˆ la densitŽ logarithmique, \cf\  \citer{Cob72} etc.) ne sont pas exploitables dans l'Žtude statistique de la relation \eqref{sbtsa}.
\par 

\goodbreak
\paraun{DŽmonstration} 
\paradeuxb{Minoration} 
Soit $k$ un entier assez grand. Posons $N_k:=  b^k-1 $ et dŽsignons par $\M_k$ l'ensemble des entiers $m$ n'excŽdant pas $N_k$ et dont le dŽveloppement en base $b$ ne contient que des chiffres $0$ et $b-1$. ƒtant donnŽ un paramtre $\varrho\in]0,1[$, nous dŽfinissons  une loi de probabilitŽ sur $\M_k$ par la formule 
$$\PP(m)=r_k(m):=\varrho^{s_b(m)/(b-1)}(1-\varrho)^{k-s_b(m)/(b-1)},$$
de sorte que
$$\eqalign{&\PP\big(s_b=j(b-1)\big)={k\choose j}\varrho^j(1-\varrho)^{k-j},\quad\sum_{m\in\M_k}\PP(m)=\sum_{0\leqslant j\leqslant k}{k\choose j}\varrho^j(1-\varrho)^{k-j}=1,\cr   &\EE(s_b)=\varrho(b-1)k.\cr}$$
\goodbreak\par \goodbreak
Par un rŽsultat classique relatif ˆ la loi binomiale, nous avons
$$\VV(s_b)=\sum_{m\in\M_k}r_k(m)\big\{s_b(m)-\varrho (b-1)k\}^2=\varrho(1-\varrho)k(b-1)^2,$$
de sorte que, en vertu de l'inŽgalitŽ de BienaymŽ-TchŽbychev,
$$\PP\Big(|s_b-\varrho (b-1)k|>T\sqrt{k}\Big) \leqslant {\varrho(1-\varrho)(b-1)^2\over T^2}\qquad (T\geqslant 1).\eqdef{BT}$$
\par 
Nous allons montrer que le dŽveloppement $a$-adique des ŽlŽments de $\M_k$ est simplement normal\note{Autrement dit, chaque chiffre y appara"t avec frŽquence asymptotique $1/a$.} sauf peut-tre pour ceux d'un sous-ensemble $\E_k$ de probabilitŽ tendant vers 0 lorsque $k\to\infty$. Cela impliquera
$$s_a(m)\sim {a-1\over 2\log a}\log N_k\qquad (m\in\M_k\sset\E_k,\,k\to\infty).$$
  Lorsque $0<\tau< 2\tau_0$, nous choisirons $\varrho=\varrho_1 $ et nous dŽduirons donc de \eqref{BT}  la relation cherchŽe $$\PP(s_b\sim \tau s_a)=1+o(1)\qquad (k\to\infty).\eqdef{s2sims3}$$
 Lorsque $\tau>\dm\tau_0$, le rŽsultat annoncŽ sera obtenu en intervertissant les r™les de $a$ et $b$. 
\par 
Notons traditionnellement $\e(u):=\e^{2\pi i u}$ $(u\in\r)$. Pour Žtudier le dŽveloppement \hbox{$a$-adique} des entiers de $\M_k$, nous introduisons les sommes de Weyl
$$\sigma_h(m,n):={1\over n}\sum_{1\leqslant\nu\leqslant n}\e\Big({hm\over a^\nu}\Big)\qquad (m\in \M_k,\,n\geqslant 1,\,h\in\z)$$ et la quantitŽ
$$\Delta_n(m):={1\over H+1}+\sum_{1\leqslant h\leqslant H}{|\sigma_h(m,n)|\over h}\cdot\eqdef{defDn}$$
Notant $\chv x$ la partie fractionnaire d'un nombre rŽel $x$, il rŽsulte alors d'une forme actuelle de la majoration classique d'Erd\H os-Tur‡n pour la discrŽpance modulo 1 d'une suite rŽelle (\cf, par exemple, \citer{Te15}, \theoreme I.6.15) que, pour tout intervalle $I\subset [0,1]$, nous avons
$$\Vbs{{1\over n}\sum_{\di{1\leqslant \nu\leqslant n}{\chv{m/a^\nu}\in I}}1-|I|}\leqslant  \Delta_n(m)\qquad (m\in\M_k,\,n\geqslant 1).$$
\par 
Notons $m=\sum_{\nu\geqslant 0}e_\nu(m)a^\nu$ le dŽveloppement $a$-adique d'un entier $m$. Ainsi
$$\chv{m/a^\nu}\in[j/a, (j+1)/a[ \,\Leftrightarrow e_{\nu-1}(m)= j\qquad (1\leqslant \nu\leqslant n,\, m\geqslant 1,\, 0\leqslant j\leqslant a-1).$$
\par 
Pour le choix $n:=\fl{(\log N_k)/\log a}$, nous avons donc
$$s_{a}(m)=\dm(a-1) n+O\big(n\Delta_n(m)\big)={(a-1)\over 2\log a}\log N_k+O\big(n\Delta_n(m)\big).\eqdef{evalsa}$$
\par
Pour complŽter la dŽmonstration, il reste donc ˆ montrer que $\Delta_n(m)$ tend vers 0 lorsque $k\to\infty$ pour presque tous les ŽlŽments $m$ de $\M_k$. Ë cette fin, nous considŽrons l'espŽrance
$$\EE\big(|\sigma_h(\cdot,n)|^2\big)=\sum_{m\in\M_k}r_k(m)|\sigma_h(m,n)|^2={1\over n^2}\sum_{1\leqslant\mu,\nu\leqslant n}\sum_{m\in\M_k}r_k(m)\e\Big(hm\Big({1\over a^\nu}-{1\over a^\mu}\Big)\Big).$$ 
La somme intŽrieure vaut
$$V_h(\mu,\nu):=\prod_{j\leqslant k}\Big\{1-\varrho+\varrho\,\e\Big(h(b-1)b^j\Big({1\over a^\nu}-{1\over a^\mu}\Big)\Big)\Big\}.$$
La minoration $|\sin \pi u|\geqslant 2\no u$ $(u\in\r)$ fournit alors
$$|V_h(\mu,\nu)|= \prod_{j\leqslant k}\Big\{1-4\varrho(1-\varrho)\sin^2\Big(\pi h(b-1)b^j\Big({1\over a^\nu}-{1\over a^\mu}\Big)\Big)\Big\}^{1/2}\leqslant \e^{- {8\varrho(1-\varrho)}S_h(\mu,\nu)},$$ 
 o l'on a posŽ $$S_h(\mu,\nu):=\sum_{0\leqslant j\leqslant k}\NO{h(b-1)b^j\Big({1\over a^\nu}-{1\over a^\mu}\Big)}^2.$$
\par 
Il rŽsulte de \eqref{defDn} que, pour tout $H\geqslant 1$,
$$\eqalign{\EE(\Delta_n)&\leqslant {1\over H+1}+\sum_{1\leqslant h\leqslant H}{1\over h}\EE\big(|\sigma_h(\cdot,n)|\big)\leqslant {1\over H}+\sum_{1\leqslant h\leqslant H}{1\over h}M_h(n),\cr}\eqdef{majEDelta}$$
o $M_h(n)\geqslant 0$ est dŽfini par $$M_h(n)^{2}={1\over n^2}\sum_{\mu,\nu\leqslant n}\e^{- 8\varrho(1-\varrho) S_h(\mu,\nu)}.\eqdef{Mhn}$$
\par 
DŽsignons par $M_h^*(n)$ la sous-somme de \eqref{Mhn} correspondant ˆ la condition supplŽ\-men\-taire $\nu<\min(  n-\sqrt{n},\mu)$. Nous avons clairement
$$M_h(n)^{2}=2M_h^*(n)+O\Big({1\over n}\Big).$$\par 
Posant $$\vartheta:={\log  a\over \log  b},\quad\varphi:={1\over \log b}, \quad J_\nu=J=\fl{\nu\vartheta}, \quad L=\fl{\varphi\log h}, \qquad \delta=\mu-\nu,$$ 
et notant que $k-J\geqslant \sqrt{k\vartheta}+O(1)$, nous pouvons donc Žcrire
 $$S_h(\mu,\nu)=
\sum_{0\leqslant  j\leqslant  k} \NO{  (b-1)b^{j-J+L+\chv{ \varphi\log h} - \chv{ \nu\vartheta}}\Big({1 }-{1\over a^{\mu-\nu}}\Big)}^2\geqslant  S_h^*(\mu,\nu)$$
avec
$$S_h^*(\mu,\nu):=\sum_{L  <  \ell\leqslant  \kappa +L } \no{b^\ell\alpha}^2,\quad \alpha=\alpha(h,\mu,\nu):=(b-1)b^{ \chv{  \varphi\log h}-\chv{ \nu\vartheta}}\big(1-{1/a^{\delta }}\big),$$
o $\kappa$ est un paramtre entier ˆ notre disposition sous la contrainte $\kappa \leqslant \dm\sqrt{k\vartheta}$. \par 
 Dans un premier temps, supposons $b\geqslant  3$.  ConsidŽrons alors le dŽveloppement $b$-adique de~$\alpha$, soit 
$$\alpha=\sum_{m \geqslant -1 }\varepsilon_m(\alpha)/b^{m}\qquad \Big(\{\varepsilon_m(\alpha)\}_{m\in\z}\in\{0,1,\ldots, b-1\}^\z\Big).$$
Notons  $\sigma_{L,\kappa }(\alpha)$ le  nombre des indices $\ell$ de $[L+1,\kappa +L ]$ tels que $\varepsilon_\ell(\alpha)=1$. 
Si $\ell$ est comptŽ dans $\sigma_{L,\kappa}(\alpha)$, nous avons
 $${1\over b}\leqslant  \chv{b^{\ell-1}\alpha}={1\over b}+\sum_{m\geqslant \ell+1}{\varepsilon_{m}(\alpha)\over b^{m-\ell+1}}\leqslant {1\over b}+\sum_{m\geqslant \ell+1}{b-1\over b^{m-\ell+1}}= {2\over b},
$$
et donc $\no{b^{\ell-1}\alpha}\geqslant 1/b$  ds que $b\geqslant 3$. Il s'ensuit que  $  S_h^*(\mu,\nu) \geqslant \sigma_{L,\kappa }(\alpha)/b^2$.
\par 
\goodbreak

Soit $s\in[0,\kappa]$. L'ensemble  $\B$  des nombres rŽels $\beta$ de $[0,(b-1)b]$ tels que $\sigma_{L,\kappa }(\beta)\leqslant  s$ est une rŽunion d'au plus  $b^{L+ 2}\sum_{0\leqslant j\leqslant  s} { \kappa \choose j}(b-1)^{\kappa -j}$ intervalles de longueurs $ \leqslant   b^{-\kappa -L}$. Pour tout $v\in]0,1]$, nous avons
$$\sum_{0\leqslant j\leqslant  s} { \kappa \choose j} (b-1)^{\kappa -j }\leqslant \sum_{0\leqslant j\leqslant \kappa}{ \kappa \choose j}(b-1)^{\kappa-j}v^{j-s}=(b-1+v)^\kappa v^{-s}.$$
Pour le choix $s=\kappa/2b $, $v= (b-1)/(2b-1) $, la majoration prŽcŽdente n'excde pas $b^{ \kappa}\e^{-\kappa/ 7b}$.
 Comme $b^L\leq h$, il s'ensuit que $\B$ est la rŽunion d'au plus  $O( hb^{ \kappa}\e^{-\kappa/ 7b})$ intervalles. 

Ainsi, nous avons
$$ S_h^*(\mu,\nu)\geqslant  {\kappa\over 2b^3}$$ 
sauf peut-tre lorsque $b^{-\chv{\nu\vartheta}} $ appartient ˆ un sous-ensemble de $[0,1]$  de mesure \hbox{$\ll \e^{-\kappa/ 7b}$} et  constituŽ d'une rŽunion d'au plus $\ll h b^\kappa\e^{-\kappa/7b}$ intervalles.   Comme $b^{-\chv{\nu\vartheta}}\geqslant 1/b$, l'ensemble exceptionnel peut-tre Žgalement caractŽrisŽ par une condition du type $\chv{\nu\vartheta}\in E$ o $E$ est un sous-ensemble de $]0,1]$ de mesure  $$|E|\ll \e^{-\kappa/ 7b}$$
et est constituŽ d'une rŽunion d'au plus $O(hb^\kappa\e^{-\kappa/ 7b} )$ intervalles.
 Le nombre des indices $\nu$ exceptionnels n'excŽdant pas $n$ est donc
$\ll \e^{-\kappa/7b }n+ hb^\kappa \e^{-\kappa/ 7b}nD_n $,
o~$D_n$ dŽsigne la discrŽpance de la suite $\{\chv{\nu\vartheta}\}_{\nu=1}^{n}$.
En sommant sur $\mu$ et $\nu$, il suit 
$$M_h(n)^{2}\ll \e^{-4\varrho(1-\varrho)\kappa/b^3}+{1\over n}+\e^{-\kappa/7b}+hb^\kappa \e^{-\kappa/ 7b}D_n. \eqdef{majMnh}$$
\par 
 Pour traiter le cas $b=2$, nous supposons $\kappa$ pair et remplaons  $\sigma_{L,\kappa }(\alpha)$  par   le  nombre $\sigma^*_{L,\kappa}(\alpha)$ des indices $j$ de l'intervalle  $  [\dm (L+1),\dm (L+ \kappa-1 )]$ tels que $\varepsilon_{2j}(\alpha) +\varepsilon_{2j+1}(\alpha)=1$. 
Pour de tels $j$, nous avons
$$\chv{2^{2j-1}\alpha}\in [\ft14, \ft34],$$ de sorte que $\no{2^{2j-1}\alpha}\geqslant 1/4$ et donc $S_h^*(\mu,\nu)\geqslant\sigma_{L,\kappa}^*(\alpha)/16$. L'ensemble $\B^*$ des nombres rŽels $\beta$ de $[0,2]$ tels que $\sigma^*_{L,\kappa}(\alpha)\leqslant s:=\kappa/ 8$ est une rŽunion d'au plus $O(hZ)$ intervalles de longueurs $\ll 1/(h2^\kappa)$ avec 
$$ Z:=2^{\kappa/2}\sum_{j\leqslant s}{\kappa/2\choose j}\leqslant 2^{\kappa/2}\sum_{0\leqslant j\leqslant \kappa/2}{\kappa/2\choose j}\big(\ft13\big)^{j-\kappa/8}=8^{\kappa/2}3^{-3\kappa/8}\leqslant 2^{11\kappa/12}. $$
En raisonnant comme prŽcŽdemment, nous obtenons donc, pour $b=2$,
$$M_h(n)^{2}\ll \e^{-\varrho(1-\varrho)\kappa/16 }+1/n+ 2^{-\kappa/12 }+h2^{ 11\kappa/12}D_n.\eqdef{majMnh-2}$$ 
\par 
ConsidŽrons la minoration \eqref{defdelta}. D'aprs la formule (8) de \citer{BT15}, nous avons
$$D_n\ll {1\over q}+{q\over n}\qquad \big(n\geqslant 
1,\,|\vartheta-r/q|\leqslant  2/q^2\big). $$
Or, d'aprs le thŽorme de Dirichlet,  il existe, pour tout $R\in[1,n]$, des entiers premiers entre eux $r$ et $q$ tels que $1\leqslant q\leqslant n/R$ et
$$\Big|\vartheta-{r\over q}\Big|\leqslant  {R\over qn}\cdot$$
En choisissant $R=n^{1/\gamma}$ et $n$ suffisamment grand, nous obtenons $$q\gg (n/R)^{1/(\gamma-1)}\geqslant  n^{1/\gamma}$$ et donc
$D_n\ll n^{-1/\gamma}.$
\par 
\goodbreak
Reportons  dans \eqref{majMnh} et \eqref{majMnh-2} en choisissant $\kappa:= 2\fl{(\log n)/( 2\gamma\log b)}$. Il vient finalement, dans tous les cas, 
$$M_h(n)\ll \sqrt{h}n^{- {3\sigma_1/2}\gamma}\ll  \sqrt{h}k^{-3\sigma_1/2\gamma},$$
o $$ \sigma_1=\sigma_1(a,b;\varrho):= {\varrho(1-\varrho)  \over 6b^3\log b }>0. $$
\par 
Sous l'hypothse $\tau<2\tau_0$, choisissons $\rho=\varrho_1=\tau/(2\tau_0)\in ]0,1[$ et  $H:=k^{ \sigma_1/  \gamma}$.
 En insŽrant les Žvaluations prŽcŽdentes dans \eqref{majEDelta}, il vient 
 $ \EE(\Delta_n)\ll k^{- \sigma_1/\gamma}$, et donc, par \eqref{evalsa}, 
 $$\PP\Big(|s_a-(b-1)k/2\tau_0|>k^{1-\sigma/\gamma}\Big)\ll  k^{-(\sigma_1-\sigma)/\gamma} ,$$
pour tout  $0<\sigma<\sigma_1(a,b)$.
Compte tenu de \eqref{BT}, il s'ensuit que \eqref{s2sims3} a bien lieu avec un terme d'erreur $\ll k^{-(\sigma_1-\sigma)/\gamma} .$
En intervertissant les r™les de $a$ et $b$ et en choisissant ˆ prŽsent $\rho=\varrho_2=\tau_0/(2\tau)$, nous obtenons le mme rŽsultat sous l'hypothse $\tau>\dm\tau_0$, quitte ˆ changer $b$ en $a$ dans la dŽfinition de $\sigma_1$.  La loi de probabilitŽ $\PP$ Žtant concentrŽe sur les entiers de $\M_k$ ayant $k\varrho$ chiffres $b-1$ dans le dŽveloppement en base $b$, nous obtenons la valeur de $c_0(a,b;\tau)$ indiquŽe dans l'ŽnoncŽ via une Žvaluation standard du coefficient binomial~${k\choose k\varrho}$.  Pour les entiers de $\M_k$ considŽrŽs, nous avons en fait
$$\min\{s_a(m),s_b(m)\}\gg k,\quad |s_a(m)-\tau s_b(m)|\ll k^{1-\sigma/\gamma},$$
d'o \eqref{majsa-tsb}. 
\medskip
\paradeuxb{Majoration}\drefdeux{maj}
Posons $\tau_a:=\dm(a-1)/\log a$, de sorte que $\tau_0:=\tau_b/\tau_a$. Nous cherchons donc ˆ majorer, lorsque $x\to\infty$, la taille $M(x)$ de l'ensemble
$$\M:=\{n\leqslant x:s_b(n)\sim \tau s_a(n)\}.$$
Supposons par exemple $\tau<\tau_0$ et introduisons un paramtre $\lambda\in]1,2[$.\par 
Soit $M^+(x):=|\{n\in\M:s_a(n)>\lambda\tau_a\log x\}|$, et $M^-(x)$ le cardinal complŽmentaire.
Nous avons, pour tout $v>0$,
$$M^+(x)\leqslant \sum_{n\leqslant x}\e^{v(s_a(n)-\lambda\tau_a\log x)}\ll x^{-v\lambda\tau_a}\Big({\e^{av}-1\over\e^v-1 }\Big)^k\ll \e^{kg(v)}$$
avec $k=(\log x)/\log a$, $g(v):=-\dm\lambda(a-1)v+\log (\e^{av}-1)-\log (\e^v-1)$. L'optimum est atteint pour
$$\dm\lambda(a-1)={a\e^{av}\over \e^{av}-1}-{\e^v\over \e^v-1}\cdot\eqdef{vlam}$$
Le membre de droite est une fonction croissante de $v$, qui vaut $\dm(a-1)$ en $v=0$ et $a-1$ ˆ l'infini. Puisque $1<\lambda<2$, l'Žquation \eqref{vlam} possde donc une solution $v=v_\lambda$. Il suit
$$M^+(x)\ll x^{g(v_\lambda)/\log a}.$$
\par 
Si $n$ est comptŽ dans $M^-(x)$, nous avons $s_b(n)\leqslant \{\lambda\tau\tau_a+o(1)\}\log x$. Posant $\mu:=\lambda\tau/\tau_0$ et supposant $\mu<1$, nous avons alors, pour tout $w>0$,
$$M^-(x)\leqslant x^{o(1)}\sum_{n\leqslant x}\e^{-w(s_b(n)-\mu\tau_b\log x)}\ll x^{w\mu\tau_b+o(1)}\Big({1-\e^{-wb}\over 1-\e^{-w}}\Big)^j\ll x^{o(1)}\e^{jh(w)}$$
avec $j=(\log x)/\log b$, $h(w):=\dm w \mu(b-1)+\log (1-\e^{-wb})-\log (1-\e^{-w}).$ L'optimum est atteint pour 
$$\dm\mu(b-1)={1\over \e^w-1}-{b\over \e^{wb}-1}\cdot\eqdef{wlam}$$
Le membre de droite est une fonction dŽcroissante de $w>0$ dont l'image est $\big]0,\dm(b-1)\big[$. Il existe donc toujours une solution $w=w_\lambda$, d'o
$$M^-(x)\ll x^{h(w_\lambda)/\log b+o(1)}.$$
Finalement, nous avons obtenu
$$M(x)\ll x^{d_0(\tau)+o(1)}$$
avec $$d_0(\tau):=\inf_{1<\lambda<\min(\tau_0/\tau,2)}\max\Big\{g(v_\lambda)/\log a,h(w_\lambda)/\log b\Big\}.$$
Comme $\dd g(v_\lambda)/\dd\lambda=-\dm(a-1)v_\lambda<0$ et $\dd h(w_\lambda)/\dd\lambda=\dm(b-1)w_\lambda>0$, la valeur $d_0(\tau)$ est atteinte lorsque $g(v_\lambda)/\log a=h(w_\lambda)/\log b$ sous rŽserve que cette Žquation possde une solution. Comme $v_1=0$, $g(0)/\log a=1$, nous avons toujours $d_0(\tau)<1$.
\par \goodbreak
ConsidŽrons le cas $a=2$, $b=3$, $\tau=1$. Nous avons alors $\tau_2=1/\log 4$, $\tau_3=1/\log 3$, $\tau_0=(\log 4)/\log 3\approx1,26186.$ La solution de \eqref{vlam} est
$$v_\lambda=\log \Big({1\over 2/\lambda-1}\Big),$$
et, notant $\Delta:=1+6\mu-3\mu^2$, la solution de \eqref{wlam} est
$$w_\lambda=\log \bigg({\sqrt{\Delta}+\mu-1\over 6}\bigg),$$
avec ˆ prŽsent $\mu=\lambda/\tau_0=\lambda(\log 3)/\log 4$.\par 
Nous avons donc
$$\eqalign {d_1(\lambda):={g(v_\lambda)\over \log 2}&=-{\lambda \over \log 4}\log \Big({1\over 2/\lambda-1}\Big)+{1\over \log 2}\log \Big({2\over 2-\lambda}\Big),\cr d_2(\lambda):={h(w_\lambda)\over \log 3}&={\lambda w_\lambda\over \log 4}+{1\over \log 3}\log \big(1-\e^{-w_\lambda}+\e^{-2w_\lambda}\big).\cr}$$
Le calcul numŽrique montre que $d_1(\lambda)=d_2(\lambda)$ pour $\lambda=\lambda_0\approx 1,0933694$. Nous avons alors $d_0(1;2,3)=d_1(\lambda_0)=d_2(\lambda_0)\approx 0.993702$.
\bigskip
\par 
\noi{\bf Remerciements.} Ce travail a bŽnŽficiŽ d'un soutien financier de l'Agence Nationale de
la Recherche (ANR) et du {\it Fonds zur Fšrderung der wissenschaftlichen
Forschung} (FWF) dans le cadre du projet bilatŽral MuDeRa
(France--Autriche), ANR-14-CE34-0009.\par 
 Le troisime auteur tient Žgalement ˆ remercier Yann Bugeaud pour d'intŽres\-santes discussions sur ce problme et pour la rŽfŽrence au travail de Cassels.

\bigskip\bigskip\goodbreak
\centerline {\twelvebf Bibliographie}
\bigskip
{\eightpoint\leftskip9truemm 

\bibtem{AS03} J.-P. Allouche \& J. Shallit, Automatic Sequences: Theory, Applications, Generalizations, Cambridge University Press, 2003.
\bibtem{B72} A. Baker, A sharpening of the bounds for linear forms in logarithms, {\it  Acta Arith.} {\bf 21} (1972), 117--129.
\bibtem{BT15} R. de la Bretche \& G. Tenenbaum, {DŽrivabilitŽ ponctuelle d'une
intŽgrale liŽe aux fonctions de Bernoulli},  
{\it  Proc. Amer. Math. Soc.} {\bf 143} (2015), no. 11, 4791--4796.
\bibtem{Ca59} J.W.S. Cassels, On a problem of Steinhaus about normal numbers,
{\it Colloq. Math. \bf7} (1959), 95Ð101.\par  
\bibtem{Cob69} A. Cobham, {On the base-dependence of sets of numbers recognizable by
finite automata}, {\it Mathematical Systems Theory} {\bf 3} (1969), 186--192.
\bibtem{Cob72} A. Cobham, {Uniform tag sequences}, {\it Math. Systems Theory} {\bf 6} (1972), 164--192.
\bibtem{DHLL17} J.-M. Deshouillers, L. Habsieger, S. Laishram, B. Landreau,  Sums of the digits in bases 2 and 3, {\it Number theoryÑDiophantine problems, uniform distribution and applications}, 211Ð217, Springer, Cham, 2017.\par 
\bibtem{S07} V. Kh. Salikhov, On the irrationality measure of $\ln 3$ (en russe), {\it Dokl. Akad. Nauk} {\bf 417} (2007), no. 6,
753--755; English translation in {\it Dokl. Math.} {\bf 76} (2007), no. 3, 955-957.
\par 
\bibtem{Te15} G. Tenenbaum, {\it Introduction ˆ la thŽorie analytique et probabiliste des
nombres}, quatrime Ždition, coll.~ƒchelles, Belin, 2015.\par
\bibtem{WW14}  Q. Wu \& L. Wang,
   On the irrationality measure of $\log3$,
   {\it J. Number Theory \bf142} (2014), 264--273.
\par 
}

\smallskip
{\leftskip-3mm\rightskip-2cm\sevenrm
\gutter=15mm \triplecolumns
 \obeylines \baselineskip=7pt
RŽgis de la Bretche
Institut de MathŽmatiques 
\hskip5mm de Jussieu-PRG
    UMR 7586
UniversitŽ Paris Diderot-Paris 7
Sorbonne Paris CitŽ, 
Case 7012, F-75013 Paris
France
\smallskip
{\eighttt regis.de-la-breteche@imj-prg.fr}
Thomas Stoll
Institut ƒlie Cartan
 UniversitŽ de Lorraine
 BP 70239
 54506 Vand{\oe}uvre-ls-Nancy Cedex
 France
 \smallskip
{\eighttt thomas.stoll@univ-lorraine.fr}
\phantom a
\phantom a
GŽrald Tenenbaum
 Institut ƒlie Cartan
 UniversitŽ de Lorraine
 BP 70239
 54506 Vand{\oe}uvre-ls-Nancy Cedex
 France
 \smallskip
 {\eighttt gerald.tenenbaum@univ-lorraine.fr}
\phantom a
\singlecolumn
\par
}

\bye